\documentclass[a4paper, 12pt]{article}
\usepackage{amssymb}

\usepackage{amsmath, amssymb, eucal, amscd, mathrsfs, txfonts}

\newtheorem{thm}{Theorem}[section]
\newtheorem{lem}[thm]{Lemma}
\newtheorem{eg}[thm]{Example}
\newtheorem{prop}[thm]{Proposition}

\newtheorem{defn}[thm]{Definition}
\newtheorem{rem}[thm]{Remark}

 \newcommand{\norm}[1]{\left\Vert#1\right\Vert}
 \newcommand{\smnoind}{{\smallskip\noindent}}
 
 \newcommand{\ti}{\tilde}

 \newcommand{\cc}{$C^*$-algebras}
 \newcommand{\ccc}{$C^*$-algebra}

 \newcommand{\s}{\textbf s}
 \newcommand{\tar}{\textbf t}
  \makeatother

\newenvironment{proof}{{\noindent \textbf{Proof:}\ }}{\hfill $\Box$\\ \smallskip}

\begin{document}

\title{$C^*$-crossed product of
groupoid actions on categories}
\author{Han Li}
\date{}
\maketitle

\begin{abstract}
Suppose that $G$ is a groupoid acting on a small category $H$ in the
sense of \cite[Definition 4]{NOT} and $H\times_\alpha G$ is the
resulting semi-direct product category (as in \cite[Proposition
8]{NOT}). We show that there exists a subcategory $H_r \subseteq H$
satisfying some nice property called ``regularity'' such that $H_r
\times_\alpha G = H\times_\alpha G$. Moreover, we show that there
exists a so-called ``quasi action'' (see Definition \ref{quasi})
$\beta$ of $G$ on $C^*(H_r)$ (where $C^*(H_r)$ is the semigroupoid
$C^*$-algebra as defined in \cite{EXE}) such that
$C^*(H_r\times_\alpha G) = C^*(H_r)\times_\beta G$ (where the
crossed product for $\beta$ is as defined in Definition
\ref{cross}).

\medskip


\smnoind Key words:  regular action, quasi action, crossed product.
\end{abstract}

\bigskip

\section{Introduction}
\bigskip

In \cite{NOT}, Ng defined the notion of an action $\alpha$ of a
small category $G$ on another small category $H$ and show that its
semi-direct product $H\times_\alpha G$ is a small category if either
$G$ is a groupoid or $G^{(0)} = H^{(0)}$ (see \cite[Propositions 8
\& 10]{NOT}). On the other hand, Exel has recently defined in
\cite{EXE} the notion of semigroupoids as well as semigroupoid
$C^*$-algebras. Since a small category is a semigroupoid, Ng asked
whether one can define a kind of ``action'' $\beta$ of $G$ on
$C^*(H)$ such that $C^*(H)\times_\beta G \cong C^*(H\times_\alpha
G)$.

\medskip

The aim of this article is to give an answer to this question. More
precisely, we will define ``quasi-actions'' of groupoids on
$C^*$-algebras and their crossed products. For any action $\alpha$
of a groupoid $G$ on a small category $H$, there exists a
subcategory $H_r\subseteq H$ such that every element $h\in H_r$ is
in the domain of some $\alpha_g$ (such an action is called
``regular'') and that $H_r\times_\alpha G = H\times_\alpha G$.
Moreover, one can define a quasi action $\ti\alpha$ of $G$ on
$C^*(H_r)$ such that the crossed product $C^*(H_r)\times_{\ti\alpha}
G$ is $*$-isomorphic to $C^*(H_r\times_\alpha G) =
C^*(H\times_\alpha G)$. On the other hand, according to
\cite[Proposition 2]{NOT}, a groupoid $G$ can be ``represented'' by
a group bundle $\{G_\xi: \xi\in G^{(0)} /R\}$ for an equivalence
relation $R$ on $G^{(0)}$. In this case, one has a decomposition
$H_r = \cup_{\xi\in G^{(0)} /R} H_\xi$ such that $\alpha$ induces a
transitive action $\ti\alpha_\xi$ of the group $G_\xi$ on $H_\xi$
for each $\xi\in G^{(0)} /R$ and that $C^*(H_r)\times_{\ti\alpha} G
\cong \bigoplus_{\xi\in G^{(0)} /R} C^*(H_\xi)\times_{\ti\alpha_\xi}
G_\xi$. Consequently, in order to understand $C^*(H\times_\alpha
G)$, one can study a collection of the crossed products of
transitive groupoid actions on the categories.

\medskip

\emph{Acknowledgements:} The author would like to thank Prof.
Chi-Keung Ng for his invaluable conversation and for pointing out
some problems in an earlier version of this work.

\bigskip

\section{\bf Preliminaries and basic definitions}
\medskip
At first, let's recall the definition of small category, whose
notations will be used through out the paper. One can find the
definition of category in \cite {CAT}.
\begin{defn}
\label{defn1}{\rm (a) A small category $\Lambda$ is a category with
its object space being a set, or equivalently, all its morphisms
form a set. We denote the object space by $\Lambda^{(0)}$, and
composable morphism pairs by $\Lambda^{(2)}$, which is a subset of
$\Lambda\times\Lambda$.\\ \\(b) Let $f, g\in\Lambda$. We shall say
that $f$ divides $g$, or that $g$ is a multiple of $f$, in symbols
$f\mid g$, if there exists $h\in\Lambda$ such that $fh=g$. We say
that $f$ and $g$ intersect if they admit a common multiple, writing
f$\Cap$g; otherwise we will say that $f$ and $g$ are disjoint,
writing $f\perp g$.}
\end{defn}If we identify a small category
with all its morphisms, every small category is a semigroupoid (see
the definition in \cite[2.1]{EXE}). It is also well known that any
groupoid is just a small category with every morphism being
invertible. Moreover, if $H$ is a small category, we set $H^{(0)}$
as its unit space, also let $\textbf s$ and $\textbf t$ be the
source and target map of each morphism being (element) respectively.
Especially, these notations are available also for the groupoid
case.

\medskip

\begin{eg}
\label{eg2}{\rm Let X be a set and R be an equivalence relation on
X. Suppose that $G_{\xi}$ is a group for any $\xi\in(X/R)$. Then (X,
R, $\{G_{\xi}\}_{\xi\in(X/R)}$) is called a group bundle over the
equivalence classes of R. Then, let $$\mathcal {G}=\{(x, g,
y):\xi\in X/R;x, y\in\xi;g\in G_{\xi}\}$$and $\textbf t, \textbf
s:\mathcal {G}\rightarrow\ X$ are defined by $\textbf t(x, g, y)=x$
and $\textbf s(x, g, y)=y$. Moreover, let
$$(x, g, y)(y, h, z)=(x, gh, z)\ and\ (x, g, y)^{-1}=(y, g^{-1}, x)$$
for any $x,y,z\in\xi$ with $xRy$ and $yRz$, and for any $g,h\in
G_\xi$. Through this way, $\mathcal {G}$ is endowed with a groupoid
structure. Indeed, every groupoid is of this type as proved in
\cite[Proposition 2]{NOT}. More precisely, for a set $X$, there is a
one to one correspondence between groupoids with unit space X and
group bundles over the equivalence classes of equivalence relations
on X, where R is defined by $x \sim y$ if $G^x_y\triangleq \textbf
t^{-1}(x)\cap \textbf s^{-1}(y)$ is non-empty for x, y in $G^{(0)}$,
and $G_{\xi}=G_{x}^{x}$ for some fixed element $x\in\xi$. The
original thought lies in \cite[1.1]{REN}.}
\end{eg}

\medskip

\medskip

\begin{eg}
\label{par}\rm Let $Par(B(\mathcal {H}))=\{(\mathcal{M}, S): S\ {\rm
partial\ isometry},\ (\ker S)^{\perp}=\mathcal{M} \}$ for some
Hilbert space $\mathcal{H}$. Define $\s(\mathcal{M}, S)=\mathcal{M}$
and $\tar(\mathcal{M}, S)=S(\mathcal{M})$. Then we have
$Par(B(\mathcal {H}))$ is a groupoid. By Gelfand-Naimark theory, any
$C^*-$algebra $A$ admits a groupoid structure on its partial
isometries denoted by $Par(A)$.
\end{eg}
\medskip

Let's recall from \cite[Definition 4]{NOT} the definition of actions
of a groupoid on a small category as well as a remark followed. The
main motivation of this definition comes from the definition of the
Lie groupoid actions on manifolds.
\begin{defn}
\label{eg4}\rm Let $G$ and $H$ be small categories. Suppose that
$\varphi: H^{(0)} \rightarrow G^{(0)}$. We let
$$G\times^\varphi H\ :=\ \{ (g,h)\in G\times H: \s(g) = \varphi(\tar(h)) = \varphi (\s(h)) \}.$$
A \emph{left action of $G$ on $H$ with respect to $\varphi$} is a
map $(g,h)\mapsto \alpha_g(h)$ from $G\times^\varphi H$ to $H$ such
that for any $(g',g)\in G^{(2)}$, $(h',h)\in H^{(2)}$ and $u\in
H^{(0)}$ with $(g,h),(g,u), (g,h')\in G\times^\varphi H$, we have:
\begin{enumerate}
\item[\rm (I).] $\alpha_g(\s(h)) = \s(\alpha_g(h))$;
\item[\rm (II).] $\alpha_g(\tar(h)) = \tar(\alpha_g(h))$;
\item[\rm (III).] $\varphi(\alpha_g(u)) = \tar(g)$;
\item[\rm (IV).] $\alpha_{\varphi(\tar(h))} (h) = h$;
\item[\rm (V).] $\alpha_{g'}(\alpha_{g}(h)) = \alpha_{g'g}(h)$;
\item[\rm (VI).] $\alpha_g(h' h) = \alpha_g(h')\alpha_g(h)$.
\end{enumerate}
For simplicity, we say that ($\varphi, \alpha$) (or just $\alpha$)
is a left action of G on H. As a convention, We will always assume
that $\varphi$ is surjective, and denote for each $g\in G$,
$$H^g\triangleq\{h\in H: \textbf s(g)=\varphi(\textbf s(h))=\varphi(\textbf
t(h))\}.$$\end{defn}

\medskip

\begin{prop}\rm\cite[Proposition 8]{NOT}
\label{rem6} Suppose that $G$ is a groupoid acting on a small
category $H$ by a left action $(\varphi, \alpha)$ and define the
semi-direct product category $H\times_\alpha G$ by
$$H\times_\alpha G\ \triangleq \ \{(h,g)\in H\times G: \tar(g) = \varphi(\s(h)) = \varphi (\tar(h))\}.$$
For any $(h,g)\in H\times_\alpha G$, we set
$$\s(h,g)\ \triangleq\ \alpha_{g^{-1}}(\s(h)) \qquad {\rm and} \qquad \tar(h,g)\ \triangleq\ \tar(h)$$
(here, we identify $u\in H^{(0)}$ with its canonical image $(u,
\varphi(u)) \in H\times_\alpha G$). Moreover, if $(h,g),(h',g')\in
H\times_\alpha G$ satisfying $\s(\alpha_{g^{-1}}(h)) = \tar(h')$, we
define
$$(h,g)(h',g')\ \triangleq\ (h\alpha_g(h'), gg').$$
This turns $H\times_\alpha G$ into a small category. If, in
addition, $H$ is a groupoid, then $H\times_\alpha G$ is also a
groupoid with
$$(h,g)^{-1} = (\alpha_{g^{-1}}(h^{-1}), g^{-1}).$$
\end{prop}

\medskip

The semi-direct product $H\times_{\alpha} G$ has a simple criterion
for the disjointness of two elements given by the following
proposition.

\begin{prop}
\label{prop7} Let G be a groupoid, and H be a small category, and
$(\varphi, \alpha)$ is a left action of G on H. For any $(h, g)$,
$(h^{'}, g^{'})\in H\times_\alpha G$, we have $(h, g)\perp$($h^{'},
g^{'}$) iff $h\perp h^{'}$ or $g\perp g^{'}$.
\end{prop}
\begin{proof}
Sufficiency is obvious, and we only prove the necessity, i. e. , if
$(h, g)\perp\ $($h^{'}, g^{'}$) and $g\Cap g^{'}$,  $h\perp h^{'}$
must hold. Otherwise, if $h\Cap h^{'}$, by definition, there exists
$k, k^{'}$ and $l, l^{'}$ such that $hk=h^{'}k^{'}$ and
$gl=g^{'}l^{'}$. Since $(h, k)\in H^{(2)}$, then $\textbf
s(h)=\textbf t(k)$, so
$$\textbf s(\alpha_{g^{-1}}(h))=\alpha_{g^{-1}}(\textbf s(h))=\alpha_{g^{-1}}(\textbf t(k))=\textbf t(\alpha_{g^{-1}}(k))$$
and
$$\textbf s(\alpha_{g^{'-1}}(h^{'}))=\alpha_{g^{, -1}}(\textbf s(h^{'}))=\alpha_{g^{'-1}}(\textbf t(k^{'}))=\textbf t(\alpha_{g^{'-1}}(k^{'})). $$
Moreover,
$$(h, g)(\alpha_{g^{-1}}(k), l)=(h\alpha_{g}\alpha_{g^{-1}}(k), gl)=(h
\alpha_{\varphi(\textbf t(k))}(k), gl)=(hk, gl). $$
$$(h^{'}, g^{'})(\alpha_{g^{'-1}}(k^{'}), l^{'})=(h^{'}\alpha_{g^{'}}\alpha_{g^{'-1}}(k^{'}), g^{'}l^{'})=(h^{'}
\alpha_{\varphi(\textbf t(k^{'}))}(k^{'}), g^{'}l^{'})=(h^{'}k^{'},
g^{'}l^{'}). $$ On the other hand, $(hk, gl)=(h^{'}k^{'},
g^{'}l^{'})$ which gives the contradiction that $(h, g)\Cap$($h^{'},
g^{'}$). This completes the proof.
\end{proof}

\medskip

Next, we define $C^*-$algebra for small category as in \cite{EXE}.
\begin{defn}\rm
\label{defn8}\cite[4.1]{EXE} Let $\Lambda$ be a small category and
let B be the unital \ccc. A mapping $S:\Lambda\rightarrow B$ will be
called a representation of
$\Lambda$ in B, if for every $f, g\in\Lambda$, \\

\begin{enumerate}
\item[\rm (I).]$S_{f}$ is a partial isometry,
\item[\rm (II).]
\[
S_{f}S_{g}=\left\{\begin{array}{ll}
 S_{fg}, &{\rm if}\ (f, g)\in\Lambda^{(2)},\\
0, &{\rm otherwise}.
\end{array}
\right.
\]
\end{enumerate}

Moreover the initial projections $Q_{f}=S_{f}^{*}S_{f}$, and the
final projections $P_{g}=S_{g}S_{g}^{*}$, are required to commute
amongst themselves and to satisfy
\begin{enumerate}
\item[\rm (III).]$P_{f}P_{g}=0, if f\perp g$,
\item[\rm (IV).]$Q_{f}P_{g}=P_{g}$, if $(f, g)\in\Lambda^{(2)}. $
\end{enumerate}
\end{defn}

\medskip

Note that we always have: $Q_{f}P_{g}=0$ if $(f,
g)\not\in\Lambda^{(2)}$. We now recall the definition of
semigroupoid \ccc\ from \cite{EXE} (see also \cite{BLA}). The
semigroupoid \ccc\ is the \ccc\ $C^{*}(\Lambda)$ generated by a
family of partial isometries $\{S_{f}\}_{f\in\Lambda}$ subject to
the relation that the correspondence $f\rightarrow S_{f}$ is a
representation with the universal property that for every
representation T of $\Lambda$ in a unital \ccc\ B there exists a
unique *-homomorphism
$$\varphi :C^{*}(\Lambda)\rightarrow B, $$
such that $\varphi(S_{f})=T_{f}$, for every $f\in\Lambda$.

\medskip

\begin{defn}\rm\label{defn10}Let G be a groupoid, and H be a small category, and $(\varphi,
\alpha)$ is a left action of G on H. We call the action {\bf
regular}, if $\varphi(\textbf s(h))=\varphi(\textbf t(h))$ for any
$h\in H$.
\end{defn}

\medskip

\begin{lem}
\label{l1}Let G, H and $(\varphi, \alpha)$ be as above. Define
$H_r\triangleq \{h\in H:\varphi(\textbf s(h))=\varphi(\textbf
t(h))\}$, then $H_r$ is a subcategory of H, and $C^*(H\times_\alpha
G)\cong C^*(H_r\times_\alpha G)$.
\end{lem}
\begin{proof}It follows from definition that $H\times_\alpha G=H_r\times_\alpha
G$.
\end{proof}

\medskip

\section{\bf Quasi $C^*-$ dynamical system and crossed product}
\medskip

Now we introduce the definition of the crossed product of a \ccc\ by
a (discrete) groupoid. \begin{defn}\rm\label{quasi}Suppose that $A$
is a \ccc\ and $G$ is a discrete groupoid. A \emph{quasi action} of
$G$ on $A$ is a map $\beta$ from $G$ to $\mathcal {N}(A)\triangleq
\{(\varphi, \mathcal{D}(\varphi))\ |\ \varphi:\ A\rightarrow A\ {\rm
is\ a\
*-homomorphism};\ \mathcal{D}(\varphi)\subseteq A\ {\rm is\ a\
closed\
*-subalgebra};\ \varphi|_{\mathcal{D}(\varphi)}:\mathcal{D}(\varphi)\rightarrow \varphi(A)\ {\rm is\ a\ *-isomorphism}\}$, satisfying

\begin{enumerate}
\item[\rm (I).] if $(s, t)\in G^{(2)}$, then
$\beta_t(A)=\mathcal{D}(\beta_s)$, $\mathcal
{D}(\beta_{st})=\mathcal {D}(\beta_t)$, and
$\beta_{st}=\beta_s\beta_t$.
\item[\rm (II).] if $(s, t)\notin G^{(2)}$, then $\beta_s\beta_t=0$.
\end{enumerate}
Also we call $(A, G, \beta)$ a $quasi\ C^{*}-dynamical\ system$.
\end{defn}

\medskip

    By definition, one always has
$\mathcal{D}(\beta_g)=\mathcal{D}(\beta_{\s(g)})$, and
$\beta_{e}|_{\mathcal{D}(\beta_e)}=id$, for any $e\in G^{(2)}$. We
now construct a quasi $C^{*}-$dynamical system by the regular action
of a groupoid on some small category.

\medskip

\begin{prop}
\label{eg9} Let $\alpha$ be a regular action of groupoid $G$ on a
small category $H$. Assume the category \ccc\ $C^{*}(H)$ is
generated by a family of partial isometries $\{S_{h}\}_{h\in H}$.\\
{\bf (a)} for every $g\in G$, the $C^*-$algebra
$C^*(H^g)=C^*(H^{\textbf s(g)})$ can be identified with the closed
*-subalgebra generated by $\{S_h\}_{h\in
H^g}$.\\
{\bf (b)} there exists a quasi action of $G$ on $C^*(H)$ such that
for any $g\in G$, $\ti{\alpha}_g|_{C^*(H^{\textbf
s(g)})}:C^*(H^{\textbf s(g)})\rightarrow C^*(H^{\textbf t(g)})$
is the *-isomorphism given by $\ti{\alpha}_g(S_h)=S_{\alpha_g(h)}$.\\
In this way, we obtain a quasi $C^*-$dynamical system $(C^*(H), G,
\ti{\alpha})$.
\end{prop}
\begin{proof}{\bf (a)} For every $g\in G^{(0)}$, let $C^*(H^g)$ be generated by
$\{T_{h}\}_{h\in H^g}$. Define $$\mathcal{L}:H^g\rightarrow C^*(H)$$
by $$\mathcal{L}(h)=S_{h}\ \ (h\in H^g).$$ One can verify that
$\mathcal{L}$ induces a representation of $H^g$ in $C^*(H)$. By the
universal property, we get a
*-homomorphism (denoted by $\ti{\mathcal{L}}$) from $C^*(H^g)$ to
$C^*(H)$ with $\ti{\mathcal{L}}(T_h)=S_h$. On the other hand, we
define
$$\mathcal{L}':H\rightarrow C^*(H^g)$$by

\[
\mathcal{L}'(h)=\left\{\begin{array}{ll}
 T_h, &t\in H^g\\
0, &{\rm otherwise}.
\end{array}
\right.
\]

Since the action is regular, we can verify that $\mathcal{L}'$
induces a representation of $H$ in $C^*(H^g)$. Again by the
universal property, we get a
*-homomorphism(denoted by $\ti{\mathcal{L}'}$) from $C^*(H)$ to $C^*(H^g)$ with
$\ti{\mathcal{L}'}(S_h)=T_h$. Obviously, $\ti{\mathcal{L}'}$ is the
left inverse of $\ti{\mathcal{L}}$, thus $\ti{\mathcal{L}}$ is isometric, hence complete the proof of {\bf (a)}.\\
\\
{\bf (b)} For every $g\in G$, define $$\ti{\alpha}_g:H\rightarrow
C^*(H^{\textbf t(g)})$$ by

\[
\ti{\alpha}_g(h)=\left\{\begin{array}{ll}
 S_{\alpha_g(h)}, &{\rm if}\ h\in H^g\\
0, &{\rm otherwise}.
\end{array}
\right.
\]
For $h_1,h_2\in H^g$, if $(h_1,h_2)\in H^{(2)}$, since
$\s(\alpha_g(h_1))=\alpha_g(\s(h_1))=\alpha_g(\tar(h_2))=\tar(\alpha_g(h_2))$,
we have $(\alpha_g(h_1), \alpha_g(h_2))\in H^{(2)}$; and if
$(h_1,h_2)\notin H^{(2)}$, we have $(\alpha_g(h_1),
\alpha_g(h_2))\notin H^{(2)}$. Moreover, if $h_1\perp h_2$, we also
have $\alpha_g(h_1)\perp \alpha_g(h_2)$. Otherwise, if
$\alpha_g(h_1)\Cap  \alpha_g(h_2)$, there exists $k_1,k_2\in G$ such
that $\alpha_g(h_1)k_1=\alpha_g(h_2)k_2$. Consider the action by
$\alpha_{g^{-1}}$ on both sides, and this gives the contradiction
that $h_1\perp h_2$. Noticing the regularity of $\alpha$ and the
fact that $\alpha_g$ keeps the composability and disjointness of any
two elements in $H^g$ , one can verify that $\ti{\alpha}_g$ is a
representation of $H$ in $C^*(H^{\textbf t(g)})$. By the universal
property, we get a
*-homomorphism (also denoted by $\ti{\alpha}_g$) from
$C^*(H)$ to $C^*(H^{\textbf t(g)})$ given by
$\ti{\alpha}_g(S_h)=S_{\alpha_g(h)}$ for $h\in H^g$. Similarly, we
can construct a *-homomorphism
$\ti{\alpha}_{g^{-1}}:C^*(H)\rightarrow C^*(H^{\textbf s(g)})$ given
by
\[
\ti{\alpha}_{g^{-1}}(S_k)=\left\{\begin{array}{ll}
 S_{\alpha_{g^{-1}}(k)}, &{\rm if}\ k\in H^{g^{-1}}\\
0, &{\rm otherwise}.
\end{array}
\right.
\]
 It is not hard to see that the maps $\ti{\alpha}_g|_{C^*(H^{\s(g)})}$ and
$\ti{\alpha}_{g^{-1}}|_{C^*(H^{\tar(g)})}$ are the inverses of each
other, so $\ti{\alpha}_g\in \mathcal{N}(C^*(H))$ holds for every
$g\in G$.  One can also verify
$$\ti{\alpha}:G\rightarrow \mathcal {N}(C^*(H))$$ sending $g\in G$ to $\ti\alpha_g$ satisfies the conditions in Definition \ref{quasi}, hence we have $\ti\alpha$
is a quasi action of $G$ on $C^*(H)$. It completes the proof of {\bf
(b)}.

\end{proof}

\medskip

\begin{defn}\rm\label{defn11}We call $(\pi, u)$ is a covariant
representation of a quasi $C^*-$dynamical system $(A,G,\beta)$, if
$\pi$ is a  *-representation of $C^*-$algebra $A$ on some Hilbert
space $\mathcal {H}$, and u is a groupoid homomorphism from G to
$Par(B(\mathcal {H}))$ with $u_su_t=0$ for any $(s,t)\notin
G^{(2)}$, satisfying the compatible conditions that
$u(g)\pi(a)u(g)^*=\pi(\beta_g(a))$ and
$u(g)\pi(a)=\pi(\beta_g(a))u(g)$ for any $a\in \mathcal
{D}(\beta_g)$.
\end{defn}

  By definition, one always has $u(g)^*=u(g^{-1})$ for any $g\in G$.

\medskip

\begin{prop}\label{cov}
Let $(A, G, \beta)$ be a quasi $C^*-$dynamical system. Assume that
$\pi$ is a *-representation on a Hilbert space $\mathcal {H}$.
Define $(\ti{\pi}, u)$ for $(A,G,\beta)$ on $l^2(G,\mathcal {H})$
by:$$(\ti\pi(a)f)(s)=\pi(\beta_{s^{-1}}(a))(f(s))\ \ \ (a\in A),$$

\[
(u_{t}f)(s)=\left\{\begin{array}{ll}
 f(t^{-1}s), &{\rm if}\ (t^{-1},s)\in G^{(2)}\\
0, &{\rm otherwise}.
\end{array}
\right.
\]
Then $(\ti\pi, u)$ is a covariant representation of $(A, G, \beta)$.
Consequently, for any quasi $C^*-$dynamical system, covariant
representations always exist.
\end{prop}
\begin{proof} It is obvious that $\ti{\pi}$ is a
*-representation of $A$ on $l^2(G,H)$. For each $t\in G$, $\ker
(u_t)=\{f\in l^2(G,H): f(s)=0\ {\rm when}\ (t, s)\in G^{(2)}\}$.
Hence $\ker(u_t)^{\perp}=\{f\in l^2(G,H): f(s)=0\ {\rm when}\ (t,
s)\notin G^{(2)}\}$, and $u_t$ acts on $\ker(u_t)^{\perp}$
isometrically, hence a partial isometry. It is not difficult to
check that $u$ induces a groupoid homomorphism, and by definition we
have $u(s)u(t)=0$ for any $(s,t)\notin G^{(2)}$. To verify
$(\ti{\pi}, u)$ is covariant representation, for any $a\in \mathcal
{D}(\beta_t)$ if $(t^{-1}, s)\in G^{(2)}$, we have
$$(u(t)\ti\pi(a)u(t^{-1})f)(s)=(\ti\pi(a)u(t^{-1})f)(t^{-1}s)=\pi(\beta_{s^{-1}t}(a))(f(s))=(\ti\pi(\beta_t(a))f)(s);$$
and if $(t^{-1}, s)\notin G^{(2)}$, we also have
$$(u(t)\ti\pi(a)u(t^{-1})f)(s)=0=\pi(\beta_{s^{-1}}\beta_t(a))(f(s))=(\ti\pi(\beta_t(a))f)(s).$$
To show that $u(t)\ti\pi(a)=\ti\pi(\beta_t(a))u(t)$, if
$(t^{-1},s)\in G^{(2)}$, we have
$$(u(t)\ti\pi(a)f)(s)=(\ti\pi(a)f)(t^{-1}s)=(\pi(\beta_{s^{-1}t}(a))(f(t^{-1}s))=(\ti\pi(\beta_t(a))u(t)f)(s);$$
and if $(t^{-1},s)\notin G^{(2)}$, we also have
$$(u(t)\ti\pi(a)f)(s)=0=(\ti\pi(\beta_t(a))u(t)f)(s).$$ It completes the proof.
\end{proof}

\medskip

 Let $(A, G, \beta)$ be as defined in Definition \ref{quasi}. (The following notations are mainly derived from
\cite[8.1]{Dav}.) In this case, the Haar system gives the counting
measure. The space of finitely supported $A-valued$ functions is the
algebra $A[G]$ of all finite sums $f=\Sigma_{t\in G}a_{t}t$ with
$a_{t}\in\mathcal {D}(\beta_{t^{-1}})$ for each $t\in G$. Whence if
$g=\Sigma_{u\in G}b_{u}u\in A[G]$, then set
$$fg=\sum_{s\in G}\left(\sum_{(t^{-1}, s)\in G^{(2)}}a_{t}\beta_{t}(b_{t^{-1}s})\right)s,
\eqno(1a)$$and$$f^{*}=\sum_{t\in G}\beta_{t}(a_{t^{-1}}^{*})t.
\eqno(1b)$$Notice that in formula($1a$), $b_{t^{-1}s}\in\mathcal
{D}(\beta_{s^{-1}t})=\mathcal {D}(\beta_{t})$, hence it makes sense.
It is not hard to prove that $a_{t}\beta_{t}(b_{t^{-1}s})\in\mathcal
{D}(\beta_{s^{-1}})$.
 Also $(s^{-1}, s)\in G^{(2)}$ and the sum always makes sense.

\medskip

\begin{prop}
\label{induce}Any covariant representation $(\pi,u)$ of quasi
$C^*-$dynamical system $(A, G, \beta)$ yields an $l^1-$contractive
*-representation of $A[G]$ by $$\sigma(f)=\sum_{t\in G}\pi(a_t)u_t.\
\ (f=\sum_{t\in G}a_t t\in A[G])$$
\end{prop}
\begin{proof}Indeed, $$\sigma(f)^*=\sum_{t\in G} u_t^*\pi(a_t)^*=\sum_{t\in G}u_{t^{-1}}\pi(a_t^*)u_tu_{t^{-1}}=\sum_{s\in G}\pi(\beta_s(a^*_{s^{-1}}))u_s=\sigma
(f^*)$$ and notice that if $(t,v)\in G^{(2)}$, then
$\s(t)=\s(v^{-1})$, so we have $b_v\in\mathcal {D}(\beta_t)$, hence
\begin{eqnarray}\sigma(f)\sigma(g)&=&\sum_{t\in G}\sum_{v\in
G}\pi(a_t)u_t\pi(b_v)u_v\nonumber\\
&=&\sum_{t\in G}\sum_{v\in G}\pi(a_t)\pi(\beta_t(b_v))u_{tv}\ \ (t,v)\in G^{(2)}\nonumber\\
&=&\sum_{s\in G}(\sum_{(t^{-1},s)\in
G^{(2)}}\pi(a_t\beta_t(b_{t^{-1}s})))u_s=\sigma(fg). \nonumber
\end{eqnarray}
Since $\pi$ is norm contractive and non-zero partial isometry has
norm 1, it follows that $\sigma$ is $l^1-$contractive. It completes
the proof.
\end{proof}

\medskip

\begin{defn}\rm\label{cross} Let $(A,G,\beta)$ be a quasi $C^*-$dynamical system.
The crossed product $A\times_\beta G$ is the enveloping \ccc\ of
$A[G]$. That is, one defines a \ccc\ norm by
$$\norm{f}=\sup_{\sigma}\norm{\sigma(f)}$$as $\sigma$ runs over all
*-representations of $A[G]$ which is $l^{1}-$contractive. Then
$A\times_{\beta}G$ is just the $C^*-$completion of $A[G]$ with this
$C^*$-norm. \end{defn}

\medskip

\begin{rem}\rm \label{cros}In Definition \ref{cross}, $\norm{f}=\sup_{\sigma}\norm{\sigma(f)}\neq
0$ if $f\in A[G]$, and $f\neq 0$. To see this, choose $\pi$ to be a
faithful *-representation of $A$ on some Hilbert space
$\mathcal{H}$, and construct a covariant representation of $(A, G,
\beta)$ as in Proposition \ref{cov}. For $g=\sum_{t\in G}a_tt\in
A[G]$ with $a_t\neq 0$, choose $h\in \mathcal{H}$ and $h\notin
\ker(\pi(\beta_{e})(a_t))$ $(e=\s(t^{-1}))$. For $f\in l^2(G,
\mathcal{H})$ defined by
\[
f(s)=\left\{\begin{array}{ll}
 h, &{\rm if}\ s=t^{-1}\\
0, &{\rm otherwise},
\end{array}
\right.
\]
one can verify that $\sigma(g)(f)(e)\neq 0$, where $\sigma$ is as
defined in Proposition \ref{induce}. Consequently, $A[G]$ is
naturally embedded into $A\times_\alpha G$.  \end{rem}

\medskip

      By Proposition \ref{induce}, Definition \ref{cross} and Remark \ref{cros}, we have that the crossed product $A\times_\beta G$ has the property that for any
covariant representation $(\pi, u)$ of a quasi $C^*-$dynamical
system $(A,G,\beta)$, there is a representation of $A\times_\beta G$
into $C^*(\pi(A), u(G))$ obtained by setting
$$\sigma(f)=\sum_{t\in G}\pi(a_t)u_t.\ \ (f=\sum_{t\in G}a_t t\in
A[G])$$

\bigskip

\section{\bf Main applications}
Let us first recall two well-known lemmas in $C^*-$algebra theory.
\begin{lem}\label{lem}(a) Let A be a $C^*-$algebra, and S be a partial
isometry in A, with final projection $SS^*=P_0$ and initial
projection $S^*S=Q_0$. Then for any projections P and Q satisfying
$S=PS=SQ$, we have $P_0\leq P$, and $Q_0\leq Q$. ($P_0$ and $Q_0$
are minimal respect to the above property) Moreover, if there is a
$T\in Par(A)$ with $TS=Q_0$, $ST=P_0$, $TT^*=Q_0$ and $T^*T=P_0$,
then  $T=S^*$.\\
\\
(b) Let $\{S_i\}_{i\in I}$ be a family of partial isometries in
$B(\mathcal{H})$ with $\{Q_i\}_{i\in I}$ and $\{P_i\}_{i\in I}$
being the initial and final projections respectively, such that for
any $i,\ j\in I$ with $i\neq j$, we have $Q_i\perp Q_j$ and
$P_i\perp P_j$ hold.  Then the sum $\sum_{i\in I}S_i$ converges
strongly, and $\sum_{i\in I}S_i\in Par(B(\mathcal{H}))$.
\end{lem}

\medskip

\begin{thm}
Let H be a small category, G a groupoid, $(\varphi, \alpha)$ is an
action of G on H. Let $H_r=\{h\in H:\varphi(\textbf
s(h))=\varphi(\textbf t(h))\}$. Then $C^{*}(H\times_{\alpha}G)\cong
C^{*}(H_r)\times_{\ti\alpha}G. $ Here $\cong$ means *-isomorphism,
and $C^*(H_r)\times_{\ti\alpha} G$ corresponds to the quasi
$C^*-$dynamical system $(C^*(H_r),G,\ti{\alpha})$ introduced in
Proposition \ref{eg9}.
\end{thm}
\begin{proof}
With the above comments, it is enough to consider the case when the
action is regular. \\
\medskip

 Let $\Lambda=H\times _{\alpha}G$. Assume that
$C^{*}(H)$ is generated by a family of partial isometries
$\{S_{h}\}_{h\in H}$, and $C^*(H\times_\alpha G)$ is generated by
$\{S_{(h,g)}\}_{(h,g)\in H\times_{\alpha} G}$. Define
$$T:\Lambda\rightarrow C^*(H)[G]\subseteq C^{*}(H)\times_{\ti\alpha}G$$ by $$T(h, g)=S_{h}g.$$
It is well defined since $(h, g)\in\Lambda$ implies that $S_h$ lies
in the domain of $\ti\alpha_{g^{-1}}$. We shall verify the
conditions in Definition \ref{defn8} to show that $T$ is actually a
representation. \\
\begin{enumerate}\item[\rm (I).]Since \begin{eqnarray}(S_{h}g)(S_{h}g)^{*}(S_{h}g)&=&(S_{h}g)(\ti\alpha_{g^{-1}}(S_{h}^{*})g^{-1})(S_{h}g)\nonumber\\
&=&S_{h}S_{h}^{*}S_{h}g=S_hg\nonumber
\end{eqnarray} we have that $T(h, g)$ is a
partial isometry.
\item[\rm (II).]If $(h_1, g_{1})(h_2, g_{2})\in \Lambda^{(2)}$, then $\s(\alpha_{g^{-1}}(h_1))=\tar(h_2)$,
and so
$(S_{h_{1}}g_{1})(S_{h_{2}}g_{2})=S_{h_{1}}S_{\alpha_{g_{1}}(h_{2})}g_{1}g_{2}$.
Suppose that $(h_1, g_{1})(h_2, g_{2})\not\in \Lambda^{(2)}$. If
$(g_{1}, g_{2})\not\in G^{(2)}$, it follows that
$(S_{h_{1}}g_{1})(S_{h_{2}}g_{2})=0$. If $(g_{1}, g_{2})\in
G^{(2)}$, we also have $(S_{h_{1}}g_{1})(S_{h_{2}}g_{2})=0$ since
$(h_1, \alpha_{g_1}(h_2))\not\in H^{(2)}$. In sum, we have
\[T(h_{1}, g_{1})T(h_{2}, g_{2})=\left\{
\begin{array}{ll}
T((h_{1}, g_{1})(h_{2}, g_{2})) & (h_{1}, g_{1})(h_{2}, g_{2})\in\Lambda^{(2)}, \\
0 & otherwise.
\end{array}
\right.
\]
Note that the initial projection for $T(h_1, g_1)$ is $Q(h_{1},
g_{1})=\ti\alpha_{g_{1}^{-1}}(Q_{h_{1}})g_{1}^{-1}g_{1}$ and the
final projection for $T(h_2, g_2)$ is $P(h_{2},
g_{2})=P_{h_{2}}g_{2}g_{2}^{-1}$. If $(g_{1}, g_{2})\in G^{(2)}$, we
have that
$g_{1}^{-1}g_{1}=g_{1}^{-1}g_{1}g_{2}g_{2}^{-1}=g_{2}g_{2}^{-1}=g_{2}g_{2}^{-1}g_{1}^{-1}g_{1}$,
hence \begin{eqnarray}Q(h_{1}, g_{1})P(h_{2},
g_{2})&=&\ti\alpha_{g^{-1}}(Q_{h_{1}})P_{h_{2}}g_{1}^{-1}g_{1}g_{2}g_{2}^{-1}\nonumber\\
&=&P_{h_{2}}\ti\alpha_{g^{-1}}(Q_{h_{1}})g_{2}g_{2}^{-1}g_{1}^{-1}g_{1}=P(h_{2}, g_{2})Q(h_{1}, g_{1}).\nonumber\\
\end{eqnarray} On the other hand, if
$(g_{1}, g_{2})\not\in G^{(2)}$, we also have $Q(h_{1},
g_{1})P(h_{2}, g_{2})=P(h_{2}, g_{2})Q(h_{1}, g_{1})=0. $\
Similarly, we can prove that initial projections $Q(h_{1}, g_{1})$
and final projections $P(h_{2}, g_{2})$ are mutually commutative.
\item[\rm (III).] When $(h_{1}, g_{1})\perp (h_{2}, g_{2})$, we
have proved in Proposition \ref{prop7} that  $h_{1}\perp h_{2}$ or
$g_{1}\perp g_{2}$, which implies that $P(h_{1}, g_{1})P(h_{2},
g_{2})=0$.
\item[\rm (IV).]Finally, if $(h_{1}, g_{1})(h_{2},
g_{2})\in\Lambda^{(2)}$, we have $(\alpha_{g_{1}^{-1}}(h_{1}),
h_{2})\in H^{(2)}$ Hence
\begin{eqnarray}Q(h_{1}, g_{1})P(h_{2},
g_{2})&=&\ti\alpha_{{g}^{-1}_{1}}(Q_{h_{1}})P_{h_{2}}g_{2}g_{2}^{-1}=Q_{\alpha_{{g}^{-1}_{1}}(h_{1})}P_{h_{2}}g_{2}g_{2}^{-1}\nonumber\\
&=&P_{h_{2}}g_{2}g_{2}^{-1}=P(h_{2}, g_{2}) \nonumber\end{eqnarray}
\end{enumerate}So we have shown that T is indeed a representation of
$\Lambda$ in $C^{*}(H)\times_{\ti\alpha}G$,and by the universal
property, there is a unique
*-homomorphism $$\ti{T}:C^{*}(H\times_{\alpha}G)\rightarrow
C^{*}(H)\times_{\ti\alpha}G$$extending T, with the property $\ti{T}(S_{(h,g)})=S_h g$.\\

\medskip

\medskip

On the other hand, we now construct a covariant representation
$(\pi, u)$ of the quasi $C^*-$dynamical system $(C^*(H), G,
\ti\alpha)$ in $B(\mathcal {H})$, where $C^*(H\times_\alpha G)$ is
represented faithfully on $\mathcal{H}$. Define $\pi:H\rightarrow
B(\mathcal {H})$, by $$\pi(h)=S_{(h, \varphi(\s(h))}.$$ We can
verify that $\pi$ defines a representation of H in $B(\mathcal
{H})$. Then by the universal property, we get a
*-homomorphism (also denoted by $\pi$) $\pi:C^*(H)\rightarrow
B(\mathcal {H})$, by $$\pi(S_h)=S_{(h, \varphi(\s(h))}.$$ Define
$u:G\rightarrow Par(B(\mathcal {H}))$, by $$u(g)=\sum_{e\in
H^{(0)}\cap H^{g^{-1}}}S_{(e, g)}.\ \ (*)$$ We now verify that
$u(g)$ is well defined. For any $g_0\in G^{(0)}$, we have
$S_{(e,g_0)}$ is a projection for any $e\in H^{g_0}$, because either
$S_{(e,g_0)}$ is 0, or otherwise an idempotent and norm 1 element.
Therefore, $S^*_{(e,g)}=S_{(\alpha_ {g^{-1}}(e), g^{-1})}$ holds. In
fact, consider $Q=S_{(\alpha_ {g^{-1}}(e), g^{-1})}S_{(e, g)}$ which
is obviously a projection, then by checking the minimality concerned
in Lemma \ref{lem} (a), that is, $S_{(e, g)}=S_{(e, g)}Q$, and for
any projection $Q'$ satisfying $S=SQ'$, we have $Q=QQ'$, hence $Q$
is the initial projection for $S_{(e, g)}$. The same procedure works
for verifying the other conditions of Lemma \ref{lem} (a). It
follows that $S^*_{(e,g)}=S_{(\alpha_ {g^{-1}}(e), g^{-1})}$. We can
now verify that the right hand side of formula (*) is strongly
convergent. Since the pairs (e, g) in the above formula are mutually
disjoint, we have $\{Q_{(e,g)}\}$ and $\{P_{(e,g)}\}$ satisfy the
conditions in lemma \ref{lem} (b). (Note that $(e_1,g)$ and
$(e_2,g)$ are disjoint if $e_1\neq e_2$, hence
$S^*_{(e_1,g)}S_{(e_2,g)}=S^*_{(e_1,g)}P_{(e_1,g)}P_{(e_2,g)}S_{(e_2,g)}=0$
by (iii) of definition \ref{defn8}, which is an important result
that we have used silently). So u(g) is well-defined. One can check
that u defines a groupoid homomorphism, and $u(g_1)u(g_2)=0$ if
$(g_1,g_2)\notin G^{(2)}$. To prove that $(\pi, u)$ is a covariant
representation of $(C^*(H), G, \ti\alpha)$, it is enough to prove
that for $a=S_h$ such that $h\in H^{g}$, we have (the following
convergence corresponds to the strong topology)
\begin{eqnarray}u(g)\pi(S_h)u(g)^*&=&S_{(\tar(\alpha_g(h)),g)}S_{(h,g^{-1}g)}S_{(\s(h),
g^{-1})}=S_{(\alpha_g(h),gg^{-1})}\nonumber\\
&=&\pi(S_{\alpha_g(h)})=\pi(\ti\alpha_g(S_h)),\nonumber\end{eqnarray}
and $u(g)\pi(S_h)=\pi(\ti\alpha_g(S_h))u(g)$ since
$$\pi(S_h)=S_{(h,g^{-1}g)}=\pi(S_h)u(g)^*u(g).$$ Notice
that $\pi(S_h)u(g)=S_{(h, gg^{-1})}S_{(\s(h),g)}=S_{(h,g)}\in
C^{*}(H\times_{\alpha}G)$ for any $S_hg\in C^*(H)[G]$, and by the
universal property, we have a
*-homomorphism
$$\ti{T{'}}:C^{*}(H)\times_{\ti\alpha}G\rightarrow
C^{*}(H\times_{\alpha}G). $$ $\ti{T{'}}$ satisfies $\ti{T{'}}(S_h
g)=S_{(h,g)}$ for any $(h,g)\in H\times_\alpha G$. It is not hard to
see that $\ti{T}$ and $\ti{T{'}}$ inverse each other, hence complete
the proof of this theorem.
\end{proof}
\medskip

Given a family of \cc\ $\{A_{\lambda}\}_{\lambda\in\Lambda}$, we
denote $\bigoplus_{\lambda}^{c_{0}}A_{\lambda}$ for those $
(a_{\lambda})\in \Pi_{\lambda\in\Lambda}A_{\lambda}$ such that for
each $\varepsilon>0$, there exists a finite subset F of $\Lambda$
for which $\norm{a_{\lambda}}<\varepsilon$ if
$\lambda\in\Lambda\setminus F$. Also
$\bigoplus_{\lambda}^{c_{00}}A_{\lambda}$ denote those
$({a_{\lambda}})$ with finite support.

\medskip

\begin{thm}
Assume that $\alpha$ is a regular action of a groupoid G on a small
category H. Let $G^{(0)}$ be the unit space of $G$, and G is
isomorphic to the groupoid given by the group bundle $(G^{(0)}, R,
\{G_{\xi}\}_{\xi\in(G^{(0)}/R)})$ over the equivalence classes of
the equivalence relations R on $G^{(0)}$, as in Example \ref{eg2}.
Let $H_{\xi}\triangleq \{h\in H:\varphi(s(h))\in G_{\xi}\}$, and
$\alpha_\xi$ be the action of $G_\xi$ on $H_\xi$ inherited from
$\alpha$. Then $C^{*}(H)\times_{\ti\alpha}G\cong
\bigoplus_{\xi}^{c_{0}}C^{*}(H_{\xi})\times_{\ti\alpha_\xi}G_{\xi}$.
\end{thm}
\begin{proof}
Assume that $C^*(H_\xi)$ is generated by a family of partial
isometries $T_h$, and $C^*(H)$ is generated by $S_h$. We define
$\pi: H_\xi \rightarrow C^*(H),$ by $$\pi(h)=S_h.$$ Then by the
universal property, we get a *-homomorphism $\ti\pi: C^*(H_\xi)
\rightarrow C^*(H),$ by $$ \pi(T_h)=S_h.$$ On the other hand, since
the action is regular, we construct a
*-homomorphism from $C^*(H)$ to $C^*(H_\xi)$ which is the left inverse of $\ti\pi$, as we did in Proposition \ref{eg9}. Thus,
we have that $\ti\pi$ is isometry, and we reasonably view
$C^*(H_\xi)$ to be a *-subalgebra of $C^*(H)$. Consider the
decomposition operator
$$\mathcal {L}:C^{*}(H)[G]\rightarrow
\bigoplus_{\xi}^{c_{00}}C^{*}(H_{\xi})[G_{\xi}].
$$Firstly, noticing that the algebraic operations on each $C^*(H_{\xi})[G_{\xi}]$ is the restriction of the operations on $C^*(H)[G]$, one can check that $\mathcal {L}$ a *-isomorphism between these two
*-algebras. Secondly, $\mathcal {L}$ is $l^1$ contractive since the $c_{00}$ norm of any element in $\bigoplus_{\xi}^{c_{00}}C^{*}(H_{\xi})[G_{\xi}]$
 is always bounded by $l^1$ norm of the corresponding element in $C^*(H) [G]$. By universal property, we get a *-homomorphism $\ti{\mathcal {L}}:C^{*}(H)\times_{\ti\alpha}G\rightarrow
 \bigoplus_{\xi}^{c_{0}}C^{*}(H_{\xi})\times_{\ti\alpha_\xi}G_{\xi}$ extending $\mathcal {L}$. Combining with density of the two *-algebras $C^{*}(H) [G]$ and $\bigoplus_{\xi}^{c_{00}}C^{*}
 (H_{\xi})[G_{\xi}]$, we have $\ti{\mathcal {L}}$ is a *-isomorphism, hence complete the proof.
\end{proof}

 Each $G_{\xi}$ above is a transitive
subgroupoid of G in the term of \cite[1.1]{REN}, that is the map (r,
d) from $G_{\xi}$ to $G_{\xi}^{(0)}\times G_{\xi}^{(0)}$ is onto;
equivalently the orbit space $G_{\xi}^{(0)}/G_{\xi}$ is single,
whence $G(x)\triangleq t^{-1}(x)\cap s^{-1}(x)$ are isomorphic for
all $x\in\xi$. It is obvious that the each action $\alpha_\xi$ is
regular. As a result, in order to study $C^*(H\times_\alpha G)$, one
can study a collection of the crossed products of regular transitive
groupoid actions on small categories.

\bigskip

\noindent School of Mathematical Sciences, Nankai University,
Tianjin 300071, China.

\smallskip\noindent
E-mail address: lihan\_math@yahoo.com

\end{document}